\documentclass[11pt]{article}
\usepackage[francais,english]{babel}
\usepackage[latin1]{inputenc}
\usepackage[top=2cm, bottom=2cm, left=2.5cm, right=2.5cm]{geometry}

\usepackage{graphicx}
\usepackage{pgf,tikz}
\usetikzlibrary{arrows}
\usepackage{amsmath}
\usepackage{enumerate}
\usepackage{float}
\usepackage{amssymb}
\usepackage{amsthm}
\usepackage{mathrsfs}
\usepackage{color}
\usepackage{fancybox}
\usepackage[retainorgcmds]{IEEEtrantools}
\usepackage{mathabx}
\usepackage{cancel}
\usepackage{listings}
\usepackage{ulem}
\usepackage{bbm}
\usepackage{epsfig}
\usepackage[colorlinks=true]{hyperref}
\usepackage{mathdots}
\newtheorem*{thm*}{Theorem}
\newtheorem*{prop*}{Proposition}
\newtheorem{defi}{Definition}

\newcommand*{\esp}{\vspace{0.5cm}}
\definecolor{xdxdff}{rgb}{0.490196078431,0.490196078431,1.}
\definecolor{wwccff}{rgb}{0.4,0.8,1.}
\definecolor{qqqqff}{rgb}{0.,0.,1.}
\definecolor{cqcqcq}{rgb}{0.752941176471,0.752941176471,0.752941176471}
\definecolor{uuuuuu}{rgb}{0.266666666667,0.266666666667,0.266666666667}
\begin{document}
\selectlanguage{english}
\title{A numerical method to solve the Stokes problem with a punctual force in source term.}
\author{Lo\"ic LACOUTURE}
\date{October 26th 2015}
\maketitle

\esp

\noindent{\bf Abstract}\vskip 0.5\baselineskip \noindent
\selectlanguage{english}
The aim of this note is to present a numerical method to solve the Stokes problem in a bounded domain with a Dirac source term, which preserves optimality for any approximation order by the finite element method. It is based on the knowledge of a fundamental solution of the associated operator over the whole space. This method is motivated by the modeling of the movement of active thin structures in a viscous fluid.


\vskip 0.5\baselineskip

\selectlanguage{francais}
\noindent{\bf R\'esum\'e} \vskip 0.5\baselineskip \noindent
{\bf Une m\'ethode num\'erique pour la r\'esolution du probl\`eme de Stokes avec une force ponctuelle en terme source.} Le but de cette note est de pr\'esenter une m\'ethode num\'erique pour la r\'esolution du probl\`eme de Stokes avec une force ponctuelle en terme source, qui assure l'optimalit\'e de l'erreur d'approximation \'el\'ements finis. Elle s'appuie sur la connaissance explicite d'une solution fondamentale de l'op\'erateur lin\'eaire associ\'e. Cette m\'ethode est motiv\'ee par la mod\'elisation du mouvement de structures fines actives dans un fluide visqueux.

Mots cl\'es : estimations d'erreur, m\'ethode \'el\'ements finis, Stokeslet, structures fines.

\selectlanguage{francais}
\section*{Version fran\c{c}aise abr\'eg\'ee}
L'\'etude du mouvement de structures fines actives dans un fluide visqueux, tels que les flagelles permettant la nage de bact\'eries ou les cils impliqu\'es dans le transport mucociliaire, conduit \`a consid\'erer le probl\`eme de Stokes avec un second membre singulier. Dans l'asymptotique d'un cil dont le diam\`etre tend vers 0 et la vitesse vers l'infini, le terme source est en fait une distribution lin\'eique de forces. Dans le but de pouvoir faire des calculs, puisque int\'egrer num\'eriquement le long d'une courbe quelconque est difficile, nous approchons la distribution lin\'eique de forces $\delta_{\Gamma}$ par une somme de forces ponctuelles $\sum c_{i} \delta_{i}$. Une preuve bas\'ee sur celle du th\'eor\`eme des sommes de Riemann permet de montrer qu'il y a convergence, au sens faible dans $H^{-3/2-s}$, pour tout $s>0$, de $\sum c_{i} \delta_{i}$ vers $\delta_{\Gamma}$ lorsque le nombre $N$ de masses de Dirac tend vers l'infini. On peut aussi pr\'eciser la convergence dans des espaces plus faibles, voir \eqref{conv}. La convergence des solutions associ\'ees se d\'eduit de l'in\'egalit\'e \eqref{convfaible}, tir\'ee de \cite{LiMa}. On est alors ramen\'e \`a l'\'etude du probl\`eme de Stokes avec une force ponctuelle en terme source.

Lorsqu'on consid\`ere un probl\`eme elliptique avec une masse de Dirac en second membre, en dimension $d \geqslant 2$, ce second membre n'\'etant pas dans $H^{-1}$, le probl\`eme sort du cadre variationnel standard bas\'e sur l'espace de Sobolev $H^1$. Si la m\'ethode des \'el\'ements finis peut \^etre d\'efinie au niveau discret, les r\'esultats de convergence classiques ne sont {\it a priori} plus valables. Dans le cas du probl\`eme de Poisson, qui peut \^etre vu comme une version scalaire et simplifi\'ee du probl\`eme de Stokes, Scott a d\'emontr\'e dans \cite{RSc} que la m\'ethode \'el\'ements finis $P_{1}$ converge en norme $\mathbb{L}^2$ \`a l'ordre 1 en 2d et 1/2 en 3d. Des estimations similaires ont \'et\'e obtenues dans \cite{HoWi} avec une m\'ethode de Galerkin discr\`ete. De plus, Apel et ses co-auteurs ont montr\'e dans \cite{ApBe} qu'en raffinant le maillage autour de la singularit\'e, on retrouvait l'ordre de convergence classique. La m\'ethode pr\'esent\'ee, bas\'ee sur la connaissance explicite d'une solution fondamentale de l'op\'erateur lin\'eaire associ\'e, fait partie d'une classe de m\'ethodes dites de {\it soustraction}, introduites en \'electroenc\'ephalographie \cite{WoKo}. Elle permet de retrouver les ordres de convergence classiques sans raffinement de maillage.

Pour fixer les id\'ees, nous allons nous int\'eresser au probl\`eme de Stokes avec des conditions aux limites de type Dirichlet homog\`enes, voir le probl\`eme \eqref{PbSto}. La particularit\'e de ce probl\`eme r\'eside en la singularit\'e du second membre : un Dirac de force appliqu\'e en un point $x_{0}$ du domaine $\Omega$. Pour cet op\'erateur, on conna\^it une solution fondamentale d\'efinie en domaine infini, appel\'ee {\it Stokeslet}, que l'on note (${\bf u}_{\delta},p_{\delta})$, voir \eqref{stokeslet}. On obtient la solution $({\bf u},p)$ du probl\`eme \eqref{PbSto} en ajoutant \`a $({\bf u}_{\delta},p_{\delta})$ un rel\`evement r\'egulier prenant ainsi en compte les conditions aux bords. La singularit\'e de la solution $({\bf u},p)$ est contenue dans la solution fondamentale $({\bf u}_{\delta},p_{\delta})$, et elle est localis\'ee au point $x_{0}$. Le principe de la m\'ethode qui suit, est de capturer cette singularit\'e pour se ramener \`a la r\'esolution d'un probl\`eme auxiliaire r\'egulier.

On commence par d\'efinir une fonction plateau $\chi$, r\'eguli\`ere, valant 1 sur un voisinage de $x_{0}$ et 0 loin de ce point, voir D\'efinition \ref{defchi}. On note ensuite ${\bf u}_{0} = \chi {\bf u}_{\delta}$ et $p_{0} = \chi p_{\delta}$, et ${\bf g}$ et $h$ les fonctions d\'efinies en \eqref{defgh}. D'apr\`es ces d\'efinitions, on remarque que les supports de ${\bf g}$ et $h$ sont contenus dans une couronne centr\'ee en $x_{0}$, voir Figure \ref{figom}. De plus, les fonctions ${\bf u}_{\delta}$ et $p_{\delta}$ \'etant analytiques en dehors de $x_{0}$, la r\'egularit\'e des fonctions ${\bf g}$ et $h$ d\'epend directement de celle de la fonction $\chi$. Finalement, pour obtenir la solution $({\bf u},p)$ de \eqref{PbSto}, il suffit de corriger les termes d'erreur ${\bf g}$ et $h$ introduits en \eqref{defgh} en r\'esolvant le probl\`eme elliptique r\'egulier \eqref{PbStoReg}, dont on note ${\bf v}$ la solution. En effet, la fonction ${\bf u}:={\bf u}_{0} + {\bf v}$ est la solution du probl\`eme \eqref{PbSto}.

Cette m\'ethode permet de passer de la r\'esolution d'un probl\`eme singulier \`a celle d'un probl\`eme auxiliaire r\'egulier. Alors que le premier converge \`a un ordre faible \cite{RSc}, le second converge \`a l'ordre optimal, quel que soit l'ordre des \'el\'ements utilis\'es. En notant ${\bf u}_{h} := {\bf u}_{0} + {\bf v}_{h}$, o\`u ${\bf v}_{h}$ est la solution num\'erique du probl\`eme \eqref{PbStoReg} obtenue par une m\'ethode \'el\'ements finis, on d\'eduit de \eqref{major} que l'erreur commise sur ${\bf u}$ est la m\^eme que celle commise sur ${\bf v}$, et on montre ainsi que la vitesse de convergence est optimale.

Par exemple, si on utilise une m\'ethode \'el\'ements finis $P_{1}$, $k=0$ suffit. On d\'efinit alors $\chi$ comme en \eqref{chiexp}, et on explicite ${\bf g}$ et $h$, valant respectivement \eqref{g2dSto} et \eqref{h2dSto} en dimension 2, et \eqref{g3dSto} et \eqref{h3dSto} en dimension~3. Apr\`es r\'esolution num\'erique du probl\`eme \eqref{PbStoReg}, on obtient finalement une solution approch\'ee ${\bf u}_{h}$ dont l'erreur $\| {\bf u} - {\bf u}_{h}\|_{\mathbb{L}^2}$ est en $O(h^2)$, quelle que soit la dimension, contre une erreur, avec une m\'ethode directe, en $O(h)$ en dimension 2 et en $O(\sqrt{h})$ en dimension 3.

Cette m\'ethode, pr\'esent\'ee dans le cas du probl\`eme de Stokes, peut se g\'en\'eraliser \`a d'autres probl\`emes elliptiques lin\'eaires, comme le probl\`eme de Poisson avec une masse de Dirac en second membre. Les conditions aux limites de type Dirichlet homog\`enes peuvent aussi \^etre remplac\'ees par des conditions de type Dirichlet non homog\`enes,  Neumann ou Robin. Enfin, la lin\'earit\'e, qui joue un r\^ole essentiel, permet en outre de r\'esoudre le cas o\`u le second membre est la somme d'un nombre fini de forces ponctuelles et d'une fonction lisse, tout en ne r\'esolvant qu'un seul probl\`eme num\'erique.

\selectlanguage{english}
% main text

\section{Introduction.}
In order to model active thin structures in a viscous fluid, such as flagella connected to bacteria or cilia involved in the mucociliary transport, we have studied the Stokes problem with a singular right-hand side. In the asymptotic of a zero diameter cilia with an infinite velocity, the source term is a lineic distribution of forces, which, in order to ease computations, will be approximated by a sum of punctual forces. After having justified this approximation, we will present a numerical method to solve the Stokes problem with a punctual force in source term, and illustrate the results by numerical simulations.

\section{Approximation of the lineic distribution of forces by a sum of punctual forces.}
Since calculating an integral on any curve is numerically very difficult, the source term, noted $\delta_{\Gamma}$, the lineic distribution of forces on a curve $\Gamma$, is approached by a sum of $N$ punctual forces $\sum c_{i} \delta_{i}$ uniformly distributed along $\Gamma$. The theorem of Riemann sums ensures that $\sum c_{i} \delta_{i}$ weakly converges to $\delta_{\Gamma}$ in $H^{-3/2-s}$, for all $s >0$. Working in weaker spaces, it is possible to adapt the proof of theorem of Riemann sums and specify the convergence :
\begin{equation}\label{conv}
\left\| \delta_{\Gamma} - \sum_{i=1}^{N} c_{i} \delta_{i} \right\|_{H^{-2-s}} \leqslant \frac{C}{\sqrt{N}} \quad \text{ and } \quad \left\| \delta_{\Gamma} - \sum_{i=1}^{N} c_{i} \delta_{i} \right\|_{H^{-5/2-s}} \leqslant \frac{C}{N}.
\end{equation}
Moreover, using a result proved by Lions and Magenes in \cite{LiMa}, which can be written in this case
\begin{equation}\label{convfaible}
\| u \|_{H^{2-r}} \leqslant C \| f \|_{H^{-r}},\ \forall r \geqslant 0,
\end{equation}
where $u$ is the solution of a regular elliptic problem with a source term $f \in H^{-r}$, we can conclude that the solution $u_{N}$ of the Stokes problem with $\sum c_{i} \delta_{i}$ right-hand side converges to the solution $u_{\Gamma}$ of the Stokes problem with $\delta_{\Gamma}$ source term, when $N$ goes to infinity. Actually, we have
\begin{equation}\label{convsol} 
\| u_{\Gamma} - u_{N} \|_{-s} \leqslant \frac{C}{\sqrt{N}} \quad \text{ and } \quad \| u_{\Gamma} - u_{N} \|_{-1/2-s} \leqslant \frac{C}{N}.
\end{equation}
Finally, the solution of the Stokes problem with a lineic distribution of forces is approached by the solution of Stokes problem with a finite sum of punctual forces in source term. By linearity and without loss of generality, in the following we will deal with a single punctual force.

\section{Numerical method to solve the Stokes problem with a Dirac source term.}\label{PartSto}
In dimension $d \geqslant 2$, the $\delta$-distribution is not continuous on $H^1$, and so the solution of an elliptic problem with Dirac source term is not regular. Consequently, classical results for the convergence of the finite element method are not valid. In the case of the Poisson problem, which can be seen as the scalar version of the Stokes problem, Scott has shown in \cite{RSc} that the $P_{1}$-finite element method converges for $\mathbb{L}^2$-norm at the order 1 in dimension 2 and at the order 1/2 in dimension 3. Similar estimates have been obtained in \cite{HoWi} with a discrete Galerkin method. Moreover, it has been shown by Apel and his co-authors \cite{ApBe} that using graded meshes, it is possible to get numerically the classical order of convergence. The aim of this section is to present a numerical method which preserves optimality for any approximation order, without using mesh grading. It is based on the knowledge of a fundamental solution of the considered linear elliptic problem. This approach fits on the class of {\it subtraction methods}, introduced in \cite{WoKo} in the context of electroencephalography.

\subsection{Principle of the method.}\label{PartStoM}
Let us consider the following problem, defined on a bounded open domain $\Omega \subset \mathbb{R}^d$,
\begin{equation}\label{PbSto}
\left\{ \begin{array}{rccl}
-\mu \triangle {\bf u} + \nabla p & = & \delta_{x_{0}} {\bf F} & \text{in } \Omega,\\
\text{div } {\bf u} & = & 0 & \text{in } \Omega,\\
{\bf u} & = & 0 & \text{on } \partial \Omega,
\end{array} \right.
\end{equation}
where $x_{0}$ is fixed in $\Omega$ and ${\bf F}$ is a vector of $\mathbb{R}^d$. Let us note that a fundamental solution of problem \eqref{PbSto} is known in dimensions 2 and 3 :
\begin{equation}\label{stokeslet}
\begin{array}{llcl}
\bullet \ d=2, \ & {\bf u}_{\delta}({\bf x}) = \frac{1}{4 \pi \mu} \left( -\ln \left( | {\bf x} | \right){\bf I}_{2} + \frac{{\bf x} {}^{t}{\bf x}}{| {\bf x} |^{2}} \right) {\bf F} & \text{ and } & p_{\delta}({\bf x}) = \frac{1}{2 \pi} \frac{{\bf x}\cdot {\bf F}}{| {\bf x} |^2},\\[3mm]
\bullet \ d=3, \ & {\bf u}_{\delta}({\bf x}) = \frac{1}{8 \pi \mu} \left( \frac{{\bf I}_{3}}{| {\bf x} |} + \frac{{\bf x} {}^{t}{\bf x}}{| {\bf x} |^{3}} \right) {\bf F}, & \text{ and } & p_{\delta}({\bf x}) = \frac{1}{4 \pi} \frac{{\bf x}\cdot {\bf F}}{| {\bf x} |^3},\\[3mm]
\end{array}
\end{equation}
where ${\bf I}_{d}$ is the identity matrix. The fundamental solution $({\bf u}_{\delta}, p_{\delta})$ does not satisfy the boundary conditions, and so it is not the solution of problem \eqref{PbSto}. But this solution can be retrieved by adding a regular lifting, therefore the whole information on the singularity of the solution $({\bf u},p)$ is contained in the fundamental solution $({\bf u}_{\delta},p_{\delta})$ and is located at $x_{0}$. In order to extract this singularity, let us fix $0<a<b<d(x_{0},\partial\Omega)$ and define $\chi$ by Definition \ref{defchi}.

\begin{figure}[H]
\begin{minipage}[c]{.55\linewidth}
\centering
\includegraphics[scale=0.5]{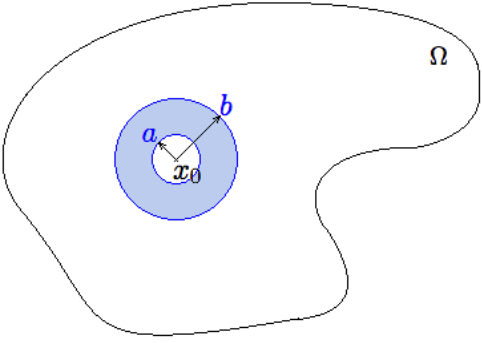}
\caption{Domain $\Omega$.}
\addtocounter{figure}{-1}
\vspace*{-0.3cm}
\caption{Domaine $\Omega$.}
\label{figom}
\end{minipage}
\begin{minipage}[c]{.40\linewidth}
\begin{defi}\label{defchi} Assume that $\chi$ is a bump function satistying for some $k \geqslant 0$,\\
\begin{itemize}
\renewcommand{\labelitemi}{$\bullet$}
\item $\chi \in H^{2+k}(\mathbb{R}^d)$,\\
\item $\chi_{\displaystyle{|_{\displaystyle{B(x_{0},a)}}}} = 1$,\\
\item $\chi_{\displaystyle{|_{\displaystyle{B(x_{0},b)^c}}}} = 0$.
\end{itemize}
\end{defi}
\end{minipage}
\end{figure}

\noindent Then, with ${\bf u}_{0} := \chi {\bf u}_{\delta}$ and $p_{0} := \chi p_{\delta}$, we define ${\bf g}$ and $h$ as
\begin{equation}\label{defgh}
{\bf g} = -\mu \triangle {\bf u}_{0} + \nabla p_{0} - \delta_{x_{0}} {\bf F} \text{ and } h = \text{div } {\bf u}_{0}.
\end{equation}
By the definitions of ${\bf u}_{\delta}$, $p_{\delta}$ and $\chi$, $supp({\bf g}) \subset \mathcal{R}^{b}_{a}(x_{0})$  and $supp(h) \subset \mathcal{R}^{b}_{a}(x_{0})$, where $\mathcal{R}^{b}_{a}(x_{0})$ is the ring centered around $x_{0}$, of internal radius $a$ and external radius $b$, see Figure \ref{figom}. Since ${\bf u}_{\delta}$ and $p_{\delta}$ are analytic on $\Omega \setminus \{ x_{0} \}$, the regularity of functions ${\bf g}$ and $h$ directly depends on the regularity of function $\chi$, namely ${\bf g} \in H^{k}(\Omega)$ and $h \in H^{k+1}(\Omega)$. Finally, it only remains to correct the terms ${\bf u}_{0}$ and $p_{0}$ by solving the regular elliptic problem
\begin{equation}\label{PbStoReg}
\left\{ \begin{array}{rccl}
-\mu \triangle {\bf v} + \nabla q & = & -{\bf g} & \text{in } \Omega,\\
\text{div } {\bf v} & = & -h & \text{in } \Omega,\\
{\bf v} & = & 0 & \text{on } \partial \Omega,
\end{array} \right.
\end{equation}
and the solution of problem \eqref{PbSto} is given by $({\bf u},p) := ({\bf u}_{0} + {\bf v} , p_{0} + q) = (\chi {\bf u}_{\delta} + {\bf v},\chi p_{\delta} + q),$ where ${\bf u}_{0}$ and $p_{0}$ are explicitly known functions and $({\bf v},q)$ is the solution of problem \eqref{PbStoReg}. Noting $({\bf v}_{h},q_{h})$ the numerical solution of problem \eqref{PbStoReg} and defining ${\bf u}_{h}:= {\bf v}_{h} + {\bf u}_{0}$ and $p_{h} = q_{h} + p_{0}$, we have,
\begin{equation}\label{major}
\begin{array}{ll}
\| {\bf u} - {\bf u}_{h} \|_{H^s(\Omega)} = \| {\bf v} - {\bf v}_{h} \|_{H^s(\Omega)}, & \text{ for } 0 \leqslant s \leqslant k +1,\\
\| p - p_{h} \|_{H^s(\Omega)} = \| q - q_{h} \|_{H^s(\Omega)}, & \text{ for } 0 \leqslant s \leqslant k .
\end{array}
\end{equation}
Actually, this method allows us to switch from the numerical computation of the solution of a singular problem with Dirac source term (with a poor convergence rate) to the numerical computation of the solution of a regular problem with an optimal convergence rate, at any required precision in terms of regularity. 

\subsection{Practical aspects.}\label{PartStoC}
For the sake of simplicity, the location of the Dirac source term will be the origin. First, we need to choose a suitable function $\chi$. Actually, to take advantage of using $P_{\ell}$-finite elements, $\ell \geqslant 1$, $\chi$ has to be $H^{\ell+1}(\mathbb{R}^d)$, in order to ensure that ${\bf g} \in H^{\ell-1}(\Omega)$ and $h \in H^{\ell}(\Omega)$, and finally to get an optimal order of convergence. For instance, for $\ell = 1$, let us define $\chi$, as a radial function, by:
\begin{equation}\label{chiexp}
\chi(r) = \left\{ \begin{array}{ll}
1 & \text{ for } r \in [0,a],\\[1mm]
\displaystyle{\frac{2r^3-3(a+b)r^2+6abr+b^2(b-3a)}{(b-a)^3}} & \text{ for } r \in [a,b],\\
0 & \text{ for } r>b,
\end{array} \right.
\end{equation}
where the function $r$ is defined on $\mathbb{R}^d$ by
\begin{equation*}
r({\bf x}) = \| {\bf x} \|_{2}.
\end{equation*}

The function ${\bf g}$ and $h$ can be explicited. According to this definition of $\chi$, ${\bf g}$ and $h$ vanish outside the ring $a < \| {\bf x} \|_{2} < b$. For $a < \| {\bf x} \|_{2} < b$, the expressions of ${\bf g}$ and $h$ depend on the dimension,
\begin{itemize}
\renewcommand{\labelitemi}{$\bullet$}
\item for $d=2$,
\begin{align} 
{\bf g}({\bf x}) & = \frac{3}{2\pi (b-a)^{3}r} \left[ \bigg( (3r^2 - 2(a+b)r + ab) \ln r + 2r^2 -2(a+b)r +2ab \bigg){\bf I}_{2} + (ab - r^2)\frac{{\bf x}{}^{t}{\bf x}}{r^2} \right] {\bf F}, \label{g2dSto} \\[3mm]
h({\bf x}) & = \frac{3(1 - \ln r)(r^2-(a+b)r+ab)}{2\pi \mu (b-a)^{3}r}  {\bf x} \cdot {\bf F}. \label{h2dSto}
\end{align}
\item for $d=3$,
\begin{align}
{\bf g}({\bf x}) & = \frac{3}{4 \pi (b-a)^3 r^2} \left[ ((a+b)r - 2r^2){\bf I}_{3} + (2ab-(a+b)r)\frac{{\bf x}{}^{t}{\bf x}}{r^2} \right] {\bf F}, \label{g3dSto}\\[3mm]
h({\bf x}) & = \frac{3(r^2 - (a+b)r+ab)}{2 \pi \mu (b-a)^3 r^2} {\bf x} \cdot {\bf F}. \label{h3dSto}
\end{align}
\end{itemize}

\section{Numerical illustrations.}
In this section, we illustrate our theorical results by a numerical example. We define $\Omega$ as the unit square and $x_{0} = (0.5, 0.5)$. The following table presents the $\mathbb{L}^2$-error for a direct method (dir. meth.) and a subtraction method (sub. meth.) respectively, for a characteristic mesh size $h$, and the estimated order of convergence (e.o.c.). Figure \ref{figerr} illustrates the section $\{ y = 0.5 \}$ of the error $|{\bf u}-{\bf u}_{h}|$ in the both cases. Numerical simulations evidence the fact that solving the auxiliary problem associated to the subtraction procedure of the singularity is more efficient than solving directly the problem with the Dirac source term.\\

\begin{center}
\begin{tabular}{  c || c  c  c  c  c  || c  }
\hline 
$h$ & $2^{-3}$ & $2^{-4}$ & $2^{-5}$ & $2^{-6}$ & $2^{-7}$ & \ e.o.c. \\ \hline
Dir. meth. & $1.02\times10^{-2}$ & $4.87\times10^{-3}$ & $2.36\times10^{-3}$ & $1.21\times10^{-3}$ & $5.89\times10^{-4}$ & 1.02 \\ 
Sub. meth. & $4.12\times10^{-3}$ & $1.33\times10^{-3}$ & $2.92\times10^{-4}$ & $6.86\times10^{-5}$ & $2.71\times10^{-5}$ & 1.88 \\ 
\hline
\end{tabular}

%\begin{figure}[h!]
%\centering
%\includegraphics[scale=0.85]{coupe.eps}
%\caption{Section $\{ x_{2} = 0\}$ of the error $|u-u_{h}|$ for the direct method and the subtraction method.}
%\label{figerr}
%\end{figure}

\begin{figure}[H]
\includegraphics[scale=0.8]{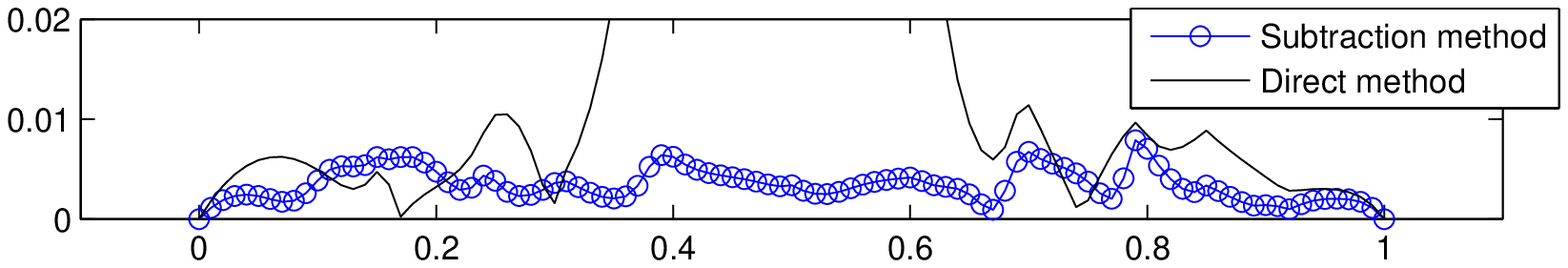}
\caption{Section $\{ y = 0.5\}$ of the error $|{\bf u}-{\bf u}_{h}|$ for the direct method and the subtraction method with $h = 0.125$.}
\addtocounter{figure}{-1}
\vspace*{-0.3cm}
\caption{Coupe $\{ y = 0.5\}$ de l'erreur $|{\bf u}-{\bf u}_{h}|$ pour la m\'ethode directe et la m\'ethode de soustraction avec $h=0.125$.}
\label{figerr}
\end{figure}
\end{center}
%
%\begin{remark}
%These examples are treated with homogeneous Dirichlet conditions. The same method is still valid in the case of non homogenous Dirichlet or any affine boundary conditions (Neumann, Robin...), up to suitable adaptations. Moreover, It is possible to generalize this method to the problem with a sum of a finite number of Dirac masses and a smooth term right-hand side.
%\end{remark}

\section{Conclusion.} 
To model active thin structures in a viscous fluid, such as flagella connected to bacteria or cilia involved in the mucociliary transport, we have studied Stokes problem with a singular right-hand side: a punctual force. However, when this problem is solved numerically, the singularity causes a poor convergence of the approximate solution to the exact solution. The method presented in this note preserves optimality for any approximation order, without using mesh grading. If the examples are treated with homogeneous Dirichlet conditions, the same method is still valid in the case of non homogenous Dirichlet or any affine boundary conditions (Neumann, Robin...), up to suitable adaptations. Similarly, the method can be generalized to the problem with a sum of a finite number of Dirac masses and a smooth term right-hand side.

\end{document}